\documentclass[11pt]{article}
\usepackage{amsfonts,amssymb,amsthm}

\newcommand{\sect}[1]{\setcounter{equation}{0}\section{#1}}
\newcommand{\subsect}[1]{\subsection{#1}}

\def\be{\begin{equation}}
\def\ee{\end{equation}}
\def\bea{\begin{eqnarray}}
\def\eea{\end{eqnarray}}

\def\C{{\mathbb C}}

\begin{document}

%%%%%%%%%%%%%%%%%%%%%%%%%%%%%%%%%%%%%%

\thispagestyle{empty}
\hfill \today

\vspace{2.5cm}

\begin{center}
\bf{\LARGE Classical Lie algebras and   Drinfeld doubles}
\end{center}

\bigskip\bigskip

\begin{center}
A. Ballesteros$^1$, E. Celeghini$^2$  and M.A. del Olmo$^3$
\end{center}

\begin{center}
$^1${\sl Departamento de F\'{\i}sica, Universidad de Burgos, \\
E-09006, Burgos, Spain.}\\
\medskip

$^2${\sl Departimento di Fisica, Universit\`a  di Firenze and
INFN--Sezione di
Firenze \\
I50019 Sesto Fiorentino,  Firenze, Italy}\\
\medskip

$^3${\sl Departamento de F\'{\i}sica Te\'orica, Universidad de
Valladolid, \\
E-47005, Valladolid, Spain.}\\
\medskip

{e-mail: angelb@ubu.es, celeghini@fi.infn.it, olmo@fta.uva.es}

\end{center}

\bigskip

\bigskip

\begin{abstract}
The Drinfeld double structure underlying the Cartan series 
$A_n, B_n, C_n, D_n$ of simple Lie algebras is discussed.
 This structure is determined by two  disjoint solvable subalgebras matched by a pairing. For the two nilpotent positive and negative root  subalgebras the pairing is natural and in the Cartan subalgebra is defined with the help of   a central extension of the algebra.
  A new   completely determined  basis     is found from    the
compatibility conditions in the double and a different perspective for quantization is presented.  Other
related Drinfeld doubles on $\C$ are also considered. 
\end{abstract}

\vskip 1cm


\vskip 0.4cm


\vfill\eject

%%%%%%%%%%%%%%%%%%%%%%%%%%%%%%%%%%%%%%%%%%%%%%%%%%%%%%%%%%%%%%% INTRODUCTION %%%%%%%%%%%%
\sect{Introduction\label{introduccion}}

The motivation of this paper is twofold.
Firstly, it is well-known that several bases are used for the Cartan subalgebra  and different authors (Bourbaki,  Chevalley,  Helgason, Weyl, ...) introduce different conventions for the commutators of root generators (see for instance Ref.~\cite{Cornwell}). Indeed, there were no  objective reasons up till now for a particular choice, we could say that no ``canonical'' basis for simple Lie algebras do not exist up to now. However,   we show here that imposing a  Drinfeld double structure on semisimple Lie algebras   the basis is completely determined. Thus, up to an overall factor,
(i) the Cartan subalgebra basis must be orthonormal and (ii) the normalization of the basis elements must be accomplished and fixed in a way different of the quoted conventions. In conclusion,  as a last (small) step of the Cartan program, a  canonical basis is found for all classical Cartan series of simple Lie algebras.  

Second, we are involved in a long term program on deformation of Lie algebras \cite{ballesteros1}. The above mentioned canonical basis, related with a coalgebra structure, constitutes the first step   for an alternative way to  quantize semisimple Lie algebras \cite{Dri87,Jimbo}. This approach allows us to consider separately the quantization of the two Borel sub-biagebras of positive and negative roots, obtaining in this way the whole quantum coalgebra, and getting afterwards the deformation.  

In Ref.~\cite{Ballesteros06} (where more technical detailes can be found) we have presented the results for the  Cartan series $A_n$. In this paper we show explicitly the canonical bases and the Drinfeld doubles for the remaining series 
$B_n, C_n, D_n$ of simple Lie algebras.  

Let us introduce the notation and basis structure of our approach.
A Lie algebra $\bar g$ is 
called a Drinfeld double \cite{Dri87} if it can be endowed  with a Manin
triple structure \cite{Dri87,majid,CP}, i.e. a set of three
Lie algebras
$(s_+,s_-,\bar g)$  such that
$s_+$ and
$s_-$ are disjoint subalgebras of $\bar g$ having the same dimension, 
$\bar g=s_+ + s_-$ as vector spaces   and  the crossed commutation rules
between $s_+$ and
$s_-$  are defined in terms of the structure tensors   of 
$s_+$ and $s_- $.  More explicitly,
let $\{ Z_p\}$ and $\{z^p\}$ be bases of  $s_+$  and $s_- $, respectively, with Lie commutators
\be\label{ces}
[Z_p,Z_q]= f^r_{p,q}\,Z_r, \qquad  [z^p,z^q]= c_r^{p,q}\,z^r   .
\ee
Provided that the following relations 
\be\label{compatibility}
c^{p,q}_r f^r_{s,t} = c^{p,r}_s f^q_{r,t}+c^{r,q}_s f^p_{r,t}
+c^{p,r}_t f^q_{s,r} +c^{r,q}_t f^p_{s,r}  \; 
\ee  
are fulfilled, a  pairing between $s_+$ and $s_-$  (i.e., a non-degenerate
symmetric
bilinear form on the vector space $s_+ + s_-$ for which $s_\pm$ are
isotropic) can be defined
\be\label{pairingz}
\langle Z_p,Z_q\rangle=0  ,\qquad
\langle Z_p,z^q\rangle=\delta_p^q , \qquad
\langle z^p,z^q\rangle=0 ,
\ee
and the remaining commutators of $\bar g$ are
\be\label{crossed}
[z^p,Z_q]= f^p_{q,r}z^r- c^{p,r}_q Z_r .
\label{zz}
\ee
 An  associated Lie bialgebra structure ($\bar g,\delta $) is also determined
by   the structure tensors of  $s_+$ and $s_-$ that
are also Lie sub-bialgebras
\be\label{delta}
\delta(Z_p)=- c_p^{q,r}\,Z_q\otimes Z_r ,\qquad 
\delta(z^p)=f^p_{q,r}\,z^q\otimes z^r .
\ee
An important property of the Drinfeld double structure is  that the cocommutator (\ref{delta})
can be derived   from  the classical
$r$-matrix $\sum_p{z^p\otimes Z_p},$
or from its skew-symmetric form 
\be\label{rmatt}
 r=\frac12 \sum_p\;{z^p\wedge Z_p}  \;.
\ee

Finally, we point out a result that turns out to be a cornerstone for our approach:  any
Drinfeld double $\bar g$ is a Lie algebra with a quadratic Casimir $C_D$
that in a certain basis $\{ Z_p,z^p\}$ can be written as
\be\label{casimirZz}
C_D=\sum_p\; [z^p,Z_p]_+ ,
\ee
where $[z^p,Z_p]_+$ denotes the anticommutator. Therefore, any (even dimensional) Lie algebra
with a quadratic Casimir in the form (\ref{casimirZz}) is a good candidate for a Drinfeld
double.

When one tries to implement the Drinfeld double structure to the Cartan series of simple Lie algebras a difficulty appears \cite{CP} for any classical algebra: a Drinfeld double can be built starting from two algebras 
$s_\pm$ isomorphic to positive and negative Borel subalgebras 
$b_\pm$. However $b_\pm$ do have the Cartan subalgebra in common and, therefore, 
cannot be identified as $s_\pm$. 

In \cite{Ballesteros06} the problem has been circumvented for $gl(n)$ by enlarging
the algebra in such a way that two disjoint solvable subalgebras
$s_\pm$, isomorphic to  Borel subalgebras, can be properly paired.
Here we follow the same approach  for all the Cartan series of
semisimple Lie algebras. Thus,
\begin{enumerate}

\item 
The $n$-dimensional Cartan subalgebra $h_n$ is 
enlarged by a central   Abelian algebra
$t_n$ generated by $I_j, (j=1,\dots,n)$. 

\item A new basis in
the $2n$-dimensional  Abelian algebra  $h_n\oplus t_n$ is defined by
\be\label{hi}
X_j:=\frac{1}{\sqrt{2}}(H_j + {\bf i} I_j), \qquad\qquad
 x^j:=\frac{1}{\sqrt{2}}(H_j- {\bf i} I_j) ,
\ee
where ${\bf i}$ is the imaginary unit.

\item The two disjoint and isomorphic solvable Lie algebras $s_+$ and
$s_-$, that contain the $X_i$'s and $x^i$'s generators and the positive
and negative roots of $g$, respectively, define a Weyl-Drinfeld double on  $\bar g=g\oplus t_n$.
\end{enumerate}

  It is worthy noting that the Drinfeld double underlying by the classical Lie algebras has the peculiarity that the structure constants 
 $c_r^{p,q}$ and 
 $f^r_{p,q}$  (\ref{ces}) verify the relation $c_r^{p,q}=-f^r_{p,q}$. So, we shall say that $\bar g$ is a Weyl-Drinfeld double (self-dual in the terminology of \cite{gomez}).

 In order to define such Weyl-Drinfeld double, we find that the usual description  of
simple Lie algebras in terms of the  Chevalley-Cartan basis (and,
obviously, Serre relations) is not suitable since, as the Killing form shows,  the Cartan subalgebra basis  is not orthonormal. On the contrary,  the Cartan subalgebra basis is   orthonormal 
 in the oscillator realization \cite{Sciarrino,Varadarajan} and for that reason we  compare our results with this last one. We shall show   that the normalization of the root vectors in the Drinfeld doubles is slightly different.
%We include once  $A_n$ because is the most relevant subalgebra in all the series.

The paper is organized as follows. 
 In section
\ref{OscillatorrealizationsofsimpleLiealgebras}, starting from  the
bosonic and  fermionic oscillator realizations of simple Lie algebras 
 we  write  the suitable
bases  for the Cartan series.
Section
\ref{DrinfelddoublesforsimpleLiealgebras}   presents the four Cartan series in terms of Drinfeld double
algebras.  A final section is devoted to comment on more general
approaches to the Drinfeld double structure of simple Lie algebras.

%%%%%%%%%%%%%%%%%%%%% SECTION II %%%%%%%%%%%%%%%%%%%%%%%%%%%%%%%%%%%%

\sect{Weyl-Drinfeld basis for classical Lie algebras}

\label{OscillatorrealizationsofsimpleLiealgebras}

%%%%%%%%%%%%%%%%
%%%%%%%%%%%%%%%%

\subsect{$A_n$ Series}

This series is the only one that supports both a bosonic as well as a fermionic oscillator
realization. So,
in terms of bosonic oscillators ($[{b_i},{b_j^\dagger}]=\delta_{ij}$)  
  the generators of $A_n$ can be written
\be\label{bosonicA}
H_i:=\frac 12 \{b_i^\dagger\, ,b_i \}   , \qquad
F_{ij}:=b_i^\dagger\,b_j  ,     \quad i\neq j  .
\ee
By using   fermionic oscillators
($\{{a_i},{a_j^\dagger}\}=\delta_{ij}$)  we have 
 \be\label{fermionicA}
 H_i:=\frac 12 [a_i^\dagger\, ,a_i ]   ,  \qquad
 F_{ij}:=a_i^\dagger\,a_j   ,  \quad i\neq j ,
\ee
where $i,j=1,\dots,n+1$. Therefore, we have $n+1$ Cartan generators $H_i$ (remember that
$\sum_i{H_i}$ is an additional central generator) and
$n(n+1)$ generators $F_{ij}$ associated to the roots.
The definitions (\ref{bosonicA}) and (\ref{fermionicA}) are slightly different from those given in \cite{Sciarrino} for a better description of the other series.

In this basis, the explicit commutation rules for $A_n$ are
\be\begin{array}{l}\label{commutatorsA}
[H_i,H_j]=0, \\[0.3cm]
[H_i,F_{jk}]=(\delta_{ij} - \delta_{ik})F_{jk}, \\[0.3cm]
[F_{ij},F_{kl}]=(\delta_{jk} F_{il}- \delta_{il} F_{kj}) + \delta_{jk} 
\delta_{il}(H_i - H_j) . 
\end{array}\ee

%%%%%%%%%%%%%%%%
%%%%%%%%%%%%%%%%

\subsect{$C_n$ Series}

The algebra  $C_n$ contains $A_{n-1}$  as a subalgebra  and its Weyl-Drinfeld basis is
given by the generators $H_i, F_{ij}$  (\ref{bosonicA}) together with
\be\label{bosonicC}
 P_{ii}:=\frac {b_i^\dagger\,b_i^\dagger}{\sqrt{2}}  ,  \quad
P_{ij}(=P_{ji}):=b_i^\dagger\,b_j^\dagger  ,  \quad
 Q_{ii}:=-\frac {b_i\,b_i} {\sqrt{2}},\quad    
Q_{ij}(=Q_{ji}):=-b_i\,b_j     ,\quad i< j ,
\ee
where $i,j=1,\dots,n$. In this way we have $n$ Cartan generators $H_i$, 
$n(n-1)$ generators $F_{ij}$, two sets of $\frac12 n(n-1)$ generators $P_{ij}$ and $Q_{ij}$ and, finally,  two sets of $n$ generators $P_{ii}$ and 
$Q_{ii}$ (where the factor $1/\sqrt{2}$, that does not appear in \cite{Sciarrino},  is imposed by the Weyl-Drinfeld invariance).

The  nonvanishing commutation rules for $C_n$ are (\ref{commutatorsA}) and
\be\begin{array}{l}\label{commutatorsC}      
[H_i,P_{jj}]= 2 \delta_{ij}\; P_{jj} ,\hskip 1.70cm
 [H_i,P_{jk}]= (\delta_{ij}+\delta_{ik})\; P_{jk} ,\quad  \\[0.3cm]  
[H_i,Q_{jj}]=- 2 \delta_{ij}\; Q_{jj} ,\hskip 1.25cm
 [H_i,Q_{jk}]= -(\delta_{ij}+\delta_{ik})\; Q_{jk}, \\[0.3cm]
 [F_{ij},P_{kk}]= \sqrt{2}  \delta_{jk}\, P_{ik}  ,\hskip1.25cm
 [F_{ij},P_{kl}]=  \delta_{jk}\, P_{il} \,
+\, \delta_{jl}\, P_{ik} +\sqrt{2} (\delta_{il}\delta_{jk}+
\delta_{ik}\delta_{jl}) P_{ii} , \\[0.3cm]
[F_{ij},Q_{kk}]= -\sqrt{2}  \delta_{ik}\, Q_{jk},\qquad
[F_{ij},Q_{kl}]=   -(\delta_{ik}\, Q_{jl} \,
+\, \delta_{il}\, Q_{jk})  -\sqrt{2} (\delta_{il}\delta_{jk}+
\delta_{ik}\delta_{jl}) Q_{jj} ,    \\[0.3cm]
[P_{ii},Q_{jj}]=2 \delta_{ij} H_i,\quad
[P_{ii},Q_{jk}]=\sqrt{2}( \delta_{ij} F_{ik}+ \delta_{ik} F_{ij}),\quad
[P_{ij},Q_{kk}]=\sqrt{2}( \delta_{ik} F_{jk}+ \delta_{jk} F_{ik}),
 \\[0.3cm]
[P_{ij},Q_{kl}]\,=\,   (\delta_{ik} F_{jl} + \delta_{jl} F_{ik}
+ \delta_{jk} F_{il} + \delta_{il} F_{jk})\,
 +  (\delta_{ik} \delta_{jl} + \delta_{jk} \delta_{il})\, (H_i+H_j).            
\end{array}\ee

%%%%%%%%%%%%%%%%
%%%%%%%%%%%%%%%%

\subsect{$D_n$ Series}

As in the preceding case, $D_n$ also contains $A_{n-1}$ as a subalgebra. However, now $D_n$ admits
a  fermionic oscillator realization. Besides the generators $H_i, F_{ij}$ 
(\ref{fermionicA}) we need two more sets of $\frac12 n(n-1)$ generators $S_{ij}$ and $T_{ij}$
given by $(i,j=1,\dots,n)$ 
\be\label{fermionicD}
S_{ij}(=-S_{ji}):=a_i^\dagger\,a_j^\dagger  ,   \qquad
T_{ij}(=-T_{ji}):=- a_i\,a_j , \qquad i<j    ,
\ee
that completes the realization of the $D_n$ algebra in the Weyl-Drinfeld basis.

The nonvanishing commutator rules for  $D_n$ are (\ref{commutatorsA}) together with  
\be\begin{array}{l}\label{commutatorsD}
[H_i,S_{jk}]\;=\; (\delta_{ij}+\delta_{ik})\; S_{jk}   , \\[0.3cm]
[H_i,T_{jk}]\;=\; -(\delta_{ij}+\delta_{ik})\; T_{jk}, \\[0.3cm]
[F_{ij},S_{kl}]\;=\;   \delta_{jk}\, S_{il} \,
-\, \delta_{jl}\, S_{ik}, \\[0.3cm]
[F_{ij},T_{kl}]\;=\; - \delta_{ik}\, T_{jl} \,
+\, \delta_{il}\, T_{jk}, \\[0.3cm]
[S_{ij},T_{kl}]\,=\, (-\delta_{jk}\, F_{il} \, - \delta_{il}\, F_{jk} \, 
+ \delta_{ik}\, F_{jl} \, + \delta_{jl}\, F_{ik}) \, + \, (\delta_{ik} \delta_{jl} - \delta_{jk}
\delta_{il})\, (H_i+H_j).
\end{array}\ee

%%%%%%%%%%%%%%%%
%%%%%%%%%%%%%%%%

\subsect{$B_n$ Series}

The Lie algebra $B_n$ contains  $D_n$ as a subalgebra. Hence, the basis of $B_n$ in the fermionic realization  will be
composed by the generators (\ref{fermionicA}), (\ref{fermionicD})
together with $2\,n$ additional generators
\be\label{fermionicB}
 U_i:=\frac{1}{\sqrt{2}}\,a_i^\dagger   ,   \qquad
 V_i:=\frac{1}{\sqrt{2}}\,a_i  , 
\ee
where, again, the factor $1/\sqrt{2}$ is required by the Weyl invariance.
The commutation rules for $B_n$ are (\ref{commutatorsA}), (\ref{commutatorsD})  plus the
non-vanishing relations involving $U$ and
$V$
\be\begin{array}{llll}\label{commutatorsB}                                                           
&[H_i,U_{j}]\;=\; \delta_{ij}\, U_j ,\qquad &
&[H_i,V_{j}]\;=\; -\delta_{ij}\, V_j, \\[0.3cm]
&[F_{ij},U_{k}]\;=\; \delta_{jk}\, U_i, \qquad &
&[F_{ij},V_{k}]\;=\; - \delta_{ik}\, V_j, \\[0.3cm]
&[T_{ij},U_{k}]\;=\;  \delta_{ik}\, V_j - \delta_{jk}\, V_i, \qquad &
&[S_{ij},V_{k}]\;=\; -  \delta_{ik}\, U_j + \delta_{jk}\, U_i, \\[0.3cm]
&[U_i,U_{j}]\;=\; \,S_{ij} , \qquad &
&[U_i,V_{j}]\;=\; \,(1 - \delta_{ij})\,F_{ij}  + \,\delta_{ij}\,H_i  , \\[0.3cm]
&[V_i,V_{j}]\;=\; \,- T_{ij} .\qquad &  &  
\end{array}\ee

%%%%%%%%%%%%%%%%%%%%%%  SECTION III %%%%%%%%%%%%%%%%%%%%%%%%%%%%%%%%%%%%%%%%%%%%%%%%%%%%%%%%%%%

\sect{Weyl-Manin triples and Lie bialgebras}

\label{DrinfelddoublesforsimpleLiealgebras}

%%%%%%%%%%%%%%%%%%%%%
%%%%%%%%%%%%%%%%%%%%%

\subsect{$A_n$ Series}

We will follow Ref.~\cite{Ballesteros06}.   Let us consider  the Lie algebra  
$gl(n+1)=A_n\oplus h$,
where $h$ is the Lie algebra generated by $\sum H_i$.

We introduce $n+1$ central generators $I_i$ and define the new generators $X_i$ and 
$ x^i$ in terms of the $H_i$ and $I_i$ as follows:
\be\label{change}
X_i:=\frac{1}{\sqrt{2}}(H_i +  {\bf i} I_i), \qquad
x^i:=\frac{1}{\sqrt{2}}(H_i- {\bf i} I_i) .
\ee

Now we consider two  solvable Lie algebras $s_+$ and $s_-$ with dimension
$(n+1)(n+2)/2$ and defined by 
\[\begin{array}{llll}
s_+:  &\{ X_i,F_{ij}\},&\qquad & i,j=1,\dots,n+1,\quad i<j,
\\[0.3cm]
s_-:  &\{ x^i,f^{ij}\},&\qquad & i,j=1,\dots,n+1,\quad i<j ,
\end{array}\]
where $f^{ij}:=F_{ji}\; ( i<j)$. Note that 
$gl(n+1)\oplus t_{n+1}= s_+ + s_-$ as vector spaces, being $t_{n+1}$ the Abelian Lie algebra generated by the $I_i$'s.
The  commutation rules for  $s_+$ and $ s_-$ are
\[\begin{array}{llll}\label{snmas}
& [X_i,X_j]=0, \qquad
& [X_i,F_{jk}]=\frac{1}{\sqrt{2}}(\delta_{ij} - \delta_{ik})\, F_{jk},\qquad
& [F_{ij},F_{kl}]=(\delta_{jk} F_{il}- \delta_{il} F_{kj}) ,\\[0.3cm]
& [x^i,x^j]=0, \qquad
& [x^i,f^{jk}]=- \frac{1}{\sqrt{2}}(\delta_{ij} - \delta_{ik})\, f^{jk},\qquad
& [f^{ij},f^{kl}]=- (\delta_{jk} f^{il}- \delta_{il} f^{kj}).
\end{array}\]

Assuming that the two algebras $s_+$ and $ s_-$ are paired by
\be\label{pairingA}
\langle x^i,X_j\rangle= \delta^i_j ,\qquad 
\langle f^{ij},F_{kl}\rangle= \delta^i_k\delta^j_l ,
\ee
we can define a bilinear form on the vector space $s_+ + s_-$ in terms of
(\ref{pairingA}) such that both $s_\pm$ are isotropic.

Taking into account (\ref{zz}) one  can easily write the crossed commutation
rules, and the compatibility relations (\ref{compatibility}) can be also checked. 
Hence, we  obtain  the Lie algebra $gl(n+1)\oplus t_{n+1}$, whose 
commutation rules in the initial basis $\{H_i,F_{ij} ,I_i \}$  are given in (\ref{commutatorsA}) plus 
$[I_i,\cdot \,]=0 $. Thus,
$(s_+,s_-, gl(n+1)\oplus t_{n+1})$ is a Weyl-Manin triple.

The canonical Lie bialgebra structure for $gl(n+1)\oplus t_{n+1}$ is determined by the 
cocommutator 
$\delta$  (\ref{delta}) and reads
\be\begin{array}{ll}\label{cocommutadorA}
\delta (I_i)=0, \\[0.3cm]
\delta (H_i)=0,\\[0.3cm]
\delta (F_{ij})=-\frac 12 F_{ij} \wedge (H_i-H_j)- \frac {{\bf i}}{2} F_{ij}
\wedge (I_i-I_j) + \sum_{k=i+1}^{j-1}{F_{ik} \wedge F_{kj}}, \qquad & i<j ,\\[0.3cm]
\delta (F_{ij})=\frac 12 F_{ij} \wedge (H_i-H_j)- \frac {{\bf i}}{2} F_{ij}
\wedge (I_i-I_j) - \sum_{k=j+1}^{i-1}{F_{ik} \wedge F_{kj}},\qquad & i>j.
\end{array}\ee
Note that  $s_+$ and its dual $s_-$ are Lie sub-bialgebras, while 
$A_n$ is a  subalgebra but not a sub-bialgebra.

The classical $r$-matrix  (\ref{rmatt})  is given by
\[ r= 
\frac12 \, \sum_{i<j}{F_{ji}\wedge F_{ij} }
+ \frac{{\bf i}}{2}\sum_{i}{H_{i}\wedge I_{i} } = r_s + r_t.
\label{rmatt3}
\]
The term $r_s$ generates the standard deformation of $gl(n+1)$ and 
$r_t$ is a twist
(not of Reshetikhin type \cite{reshetikhin}).
When all the $I_i$ are equal the twist
$r_t$ becomes trivial.

It is worth noticing  that in this construction  the chain
$gl(n)\oplus t_{n} \subset gl(n+1)\oplus t_{n+1} $ is preserved at the level of Lie bialgebras.

%%%%%%%%%%%%%%%%%%%%%%%%
%%%%%%%%%%%%%%%%%%%%%%%%

\subsect{$C_n$ Series} 

The two subalgebras $s\pm$, isomorphic to the  Borel subalgebras $b\pm$ of $C_n$, are
\[\begin{array}{llll}
s_+:  &\{ X_i,F_{ij},P_{ij},P_{ii}\},&\qquad & i,j=1,\dots,n,\quad i<j,
\\[0.3cm]
s_-:  &\{ x^i,f^{ij},p^{ij},p^{ii}\},&\qquad & i,j=1,\dots,n,\quad i<j  ,
\end{array}\]
 where the $X_i,x^i,F_{ij}$ and $f^{ij}$ were defined in the preceding section and $p^{ij}:=Q_{ij}$ and $p^{ii}:=Q_{ii}$. The Lie commutators of
the subalgebras $s_+$ and
$s_-$ are obtained from  (\ref{commutatorsA}) and (\ref {commutatorsC}).

The pairing
\[\label{pairingC}
\langle x^i,X_j\rangle= \delta^i_j ,\qquad 
\langle f^{ij},F_{kl}\rangle= \delta^i_k\delta^j_l ,
\qquad \langle p^{ij},P_{kl}\rangle= \delta^i_k\delta^j_l ,
\qquad \langle p^{ii},P_{jj}\rangle= \delta^i_j. 
\]
allows us to derive the crossed   commutators  (\ref{zz}) in this case. It can be proven that the full
set of relations are nothing but the commutation rules for $C_n\oplus t_n$.
Therefore, we can avoid to check the compatibility conditions (\ref{compatibility}) since 
$C_n\oplus t_n$ is a well known
Lie algebra, and we conclude that
$(s_+,s_-,C_n\oplus t_n)$ is a Weyl-Manin triple.

The cocommutator $\delta$  (\ref{delta}) for the generators $H_i, F_{ij}$ is given by 
(\ref{cocommutadorA}) and for the remaining generators read
\[\begin{array}{ll}\label{bialgebra gln}
\delta (P_{ii})=(H_i + {\bf i}\,I_i)\wedge P_{ii} +
\sqrt{2}\,\displaystyle{\sum_{k>i}{F_{ik}\wedge P_{ik} }},
\\[0.3cm]
\delta (P_{ij})=\frac12 [(H_i + H_j) + {\bf i}\,(I_i + I_j)]\,\wedge P_{ij} + \sqrt{2}\,F_{ij}\wedge
P_{jj}+
\displaystyle{\sum_{m>i,m\neq j}{F_{im}\wedge P_{mj}}}
,\qquad & i< j ,
\\[0.3cm]
\delta (Q_{ii})=(H_i - {\bf i}\,I_i)\wedge P_{ii} +
\sqrt{2}\,\displaystyle{\sum_{k>i}{F_{ki}\wedge Q_{ik} }},
\\[0.3cm]
\delta (Q_{ij})=\frac12 [(H_i + H_j) - {\bf i}\,(I_i + I_j)]\,\wedge Q_{ij} + \sqrt{2}\,F_{ji}\wedge
Q_{jj}+
\displaystyle{\sum_{m>i,m\neq j}{F_{mi}\wedge Q_{mj}}}.
\end{array}\]
As in the previous case, $s_\pm$ are Lie sub-bialgebras. 

The  $r$-matrix  (\ref{rmatt}),  in the basis
$\{H_i,F_{ij} ,I_i ,P_{ij},Q_{ij}\}$, 
is written as
\[ r= 
\frac12 \, \left( \sum_{i<j}{F_{ji}\wedge F_{ij} }
 + \sum_{i\leq j}{Q_{ij}\wedge P_{ij} } \right) + \frac{{\bf i}}{2}\sum_{i}{H_{i}\wedge I_{i} }=
 r_s + r_t. 
\label{rmatt3}
\]
Again $r_s$ generates the standard deformation of $C_n$ and 
$r_t$ is a twist.
When all the $I_i$ are zero 
$r_t$ vanishes. 

Note also that the chain
$C_n\oplus t_{n} \subset C_{n+1}\oplus t_{n+1}  $  is preserved at the level of Lie bialgebras. 

%%%%%%%%%%%%%%%%
%%%%%%%%%%%%%%%%%

\subsect{$D_n$ Series}

The  subalgebras $s_+$ and $s_-$ are
\[\begin{array}{llll}
s_+:  &\{ X_i,F_{ij},S_{ij}\},&\qquad & i,j=1,\dots,n,\quad i<j,
\\[0.3cm]
s_-:  &\{ x^i,f^{ij},s^{ij}\},&\qquad & i,j=1,\dots,n,\quad i<j ,
\end{array}\]
where the $X_i$, $x^i$ and $f^{ij}$ were defined in the preceeding sections and with $s^{ij}:=T_{ij}$.
The  commutators of $s_+$ and $s_-$ are obtained from  (\ref{commutatorsA}) 
and (\ref {commutatorsD}).

The crossed   commutators  (\ref{zz})
are obtained by making use of the pairing
\[\label{pairingD}
\langle x^i,X_j\rangle= \delta^i_j ,\qquad 
\langle f^{ij},F_{kl}\rangle= \delta^i_k\delta^j_l ,
\qquad \langle s^{ij},S_{kl}\rangle= \delta^i_k\delta^j_l .
\]
The complete set of Lie commutators allows us to identify  $s_+ +s_-$ with the Lie algebra  $D_n\oplus t_n$, hence
we  can avoid to check the compatibility conditions (\ref{compatibility}).  So,
$(s_+,s_-,D_n\oplus t_n)$ is a Weyl-Manin triple.

The cocommutator $\delta$  (\ref{delta}) for the generators $H_i, F_{ij}$ is given by 
(\ref{cocommutadorA}) and for the remaining generators takes the form
\be\begin{array}{ll}\label{cocommutadorD}
\delta (S_{ij})=
\frac12 [(H_i + H_j) + {\bf i}\,(I_i + I_j)]\,\wedge S_{ij} + 
\displaystyle{\sum_{k>i,k\neq j}{F_{ik}\wedge S_{kj}}},
\qquad & i< j ,
\\[0.3cm]
\delta (T_{ij})=
\frac12 [(H_i + H_j) - {\bf i}\,(I_i + I_j)]\,\wedge T_{ij} + 
\displaystyle{\sum_{k>i,k\neq j}{F_{ki}\wedge T_{kj}}},
\qquad & i< j.
\end{array}\ee
The subalgebras $s_\pm$ are Lie sub-bialgebras. 
The chain
$D_n\oplus t_{n} \subset D_{n+1}\oplus t_{n+1}  $ is again preserved at the level of Lie
bialgebras.

The classical $r$-matrix  (\ref{rmatt}),  is given by 
\[ r= 
\frac12 \left( \sum_{i<j}{F_{ji}\wedge F_{ij} }
 +  \, \sum_{i< j}{T_{ij}\wedge S_{ij} }\right) + \frac{{\bf i}}{2}\sum_{i}{H_{i}\wedge I_{i} } = r_s +
r_t ,
\label{rmatt3}
\]
where $r_s$ generates the standard deformation of $D_n$ and 
$r_t$ is a twist that vanishes when all the $I_i\to 0$.

 %%%%%%%%%%%%%%%%
 %%%%%%%%%%%%%%%%%

\subsect{$B_n$ Series}

The subalgebras $s_+$ and $s_-$ are
\[\begin{array}{llll}
s_+:  &\{ X_i,F_{ij},S_{ij},U_i\},&\qquad & i,j=1,\dots,n,\quad i<j,
\\[0.3cm]
s_-:  &\{ x^i,f^{ij},s^{ij},u^i\},&\qquad & i,j=1,\dots,n,\quad i<j ,
\end{array}\] 
where  
$s^{ij}:=T_{ij}$ and $ u^i:=V_i$.
Their commutators are obtained from (\ref{commutatorsA}) 
and (\ref{commutatorsB}).

The appropriate pairing for this case is
\[\label{pairingB}
\langle x^i,X_j\rangle= \delta^i_j ,\qquad 
\langle y^{ij},Y_{kl}\rangle= \delta^i_k\delta^j_l ,
\qquad \langle s^{ij},S_{kl}\rangle= \delta^i_k\delta^j_l ,
\qquad 
\langle u^i,U_j\rangle= \delta^i_j .\]

The complete set of Lie commutators allows us to identify  $s_+ +s_-$ with the Lie algebra  $B_n\oplus t_n$, 
so we can avoid to check the compatibility conditions (\ref{compatibility}).  Hence,
$(s_+,s_-,B_n\oplus t_n)$ is a Weyl-Manin triple.

The cocommutator $\delta$  (\ref{delta}) for the generators $H_i, F_{ij}, S_{ij},T_{ij}$ is given by 
(\ref{cocommutadorA}) and  (\ref{cocommutadorD}) and for the remaining generators reads
\[\begin{array}{ll}\label{cocommutadorB}
\delta (S_{ij})=
\frac12 [(H_i + H_j) + {\bf i}\,(I_i + I_j)]\,\wedge S_{ij} + 
\displaystyle{\sum_{k>i,k\neq j}{F_{ik}\wedge S_{kj}}} + U_i\wedge U_j ,
\qquad & i< j ,
\\[0.3cm]
\delta (T_{ij})=
\frac12 [(H_i + H_j) - {\bf i}\,(I_i + I_j)]\,\wedge T_{ij} + 
\displaystyle{\sum_{k>i,k\neq j}{F_{ki}\wedge T_{kj}} }+ V_i\wedge V_j ,
\qquad & i< j ,
\\[0.3cm]
\delta (U_{i})=
\frac12 (H_i  + {\bf i}\, I_i )\wedge U_i + \displaystyle{\sum_{i< k \leq
n}{F_{ik}\wedge U_{k}}} ,\qquad & 
\\[0.3cm]
\delta (V_{i})=\frac12 (H_i  - {\bf i}\, I_i )\wedge V_i +
\displaystyle{\sum_{1\leq k <i}{F_{ki}\wedge V_{k}}} .\qquad & 
\end{array}\]
Also in this case $s_\pm$ are Lie sub-bialgebras. 
The chain $B_n\oplus t_n\subset B_{n+1}\oplus t_{n+1}$ is preserved at the level of Lie bialgebras.
Note that, while $D_n\subset B_n$ as algebras,  $D_n\oplus t_n\not\subset B_{n}\oplus t_{n}$ as
bialgebras.

The classical $r$-matrix  (\ref{rmatt}) is written as
\[ r= 
\frac12 \left( \sum_{i<j}{F_{ji}\wedge F_{ij} }
+  \, \sum_{i< j}{T_{ij}\wedge S_{ij} } +
\sum_{i}{V_{i}\wedge U_{i} } \right) + \frac{{\bf i}}{2}\sum_{i}{H_{i}\wedge I_{i} } = r_s +
r_t.
\label{rmatt3}
\]
Again $r_s$ generates the standard deformation of $B_n$ and 
$r_t$ is the twist.

%%%%%%%%%%%%%%%%%%%%%%%%%%%%%%%%%%%%%%%%%%%%%%%%%%%%%%%%%%%%%%%%%%%%%

\sect{Other decompositions in $\C$}

The canonical Weyl-Drinfeld double structures above
introduced are not the only possible ones. 
Indeed, we can consider other possibilities by taking into account that the quadratic Casimir,
in the Cartan-Weyl basis,   can be always written as
\be\label{quadraticcasimir}
{\cal C}_2= \sum H_i^2+  \sum\; [X^+_j,X^-_j]_+ ,
\ee
where the $H_i$'s determine a well defined basis in the Cartan subalgebra and $\{X^+_j\}$ (resp. $\{X^-_j\}$) constitues an appropriate basis for  the nilpotent algebra of
positive (resp. negative) root generators.
Obviously, the roots $\{X^+_j\}$ (resp. $\{X^-_j\}$) can be immediately associated to the
subalgebra generated by the
$\{Z_p\}$ (resp. $\{z^p\}$). Therefore, the problem is to
fit the Cartan generators within the scheme. 

In the previous section we have given a  solution to this problem for any
semisimple Lie algebra, by enlarging the ($n$ dimensional) Cartan subalgebra through the
addition of
$n$ central generators  $\{I_j\}$ and by taking into account that
\be
\sum_{i}^n (H_i^2 + I_i^2) =  \sum_{i}^n \; [\frac{1}{\sqrt{2}}(H_i + {\bf i} I_i),\frac{1}{\sqrt{2}}(H_i -{\bf i} I_i)]_+.
\ee
At this point, the full
Drinfeld double structure comes out in a natural way by including the
$\frac{1}{\sqrt{2}}(H_i \pm {\bf i} I_i)$ generators in the
$s_\pm$ subalgebras 
\[
s_\pm= \{\frac{1}{\sqrt{2}} (H_i\pm {\bf i} I_{i}), X^\pm_j\}.
\]

For the Lie algebra $A_1$  this is the only possible solution.
However, for Lie algebras whose rank is even, it is possible to construct a different Drinfeld
double without introducing any additional central operator $I_i$.  For instance, in the case
of 
$A_2$, by  using the Gell-Mann basis \cite{Bacry}, the two solvable algebras can be chosen as 
\[
s_\pm=\{ \lambda_3 \pm {\bf i} \lambda_8,  \lambda_1 \pm {\bf i} \lambda_2, \lambda_4 \pm {\bf i} \lambda_5, \lambda_6 \pm {\bf i} \lambda_7\}.
\]

Another example of double structure can be constructed for the $D_2\thickapprox
A_1\oplus A_1$ algebra through the  solvable
algebras
\[
s_\pm=\{J_3\pm {\bf i} L_3,J_1\pm {\bf i} J_2,L_1\pm {\bf i} L_2\},
\]
where the generators $J_i$'s belong to the first $A_1$ algebra and the $L_i$'s to the
second one.

However, for odd dimensional algebras at least one additional generator must be introduced in order to get a global even
dimension. In general, all the intermediate cases among the  canonical case and the
previous ones can be considered by introducing an  algebra $t_m$  of central elements
$I_i$ with
$0\leq m\leq n$ and  such that $(n-m)/2 \in Z^+$. Thus,  the two  solvable algebras can be
defined as
\[
s_\pm= \{\frac{1}{\sqrt{2}} (H_i\pm {\bf i} H_{j}),\frac{1}{\sqrt{2}} (H_k\pm {\bf i} I_{k}), X^\pm_l \},
\]
where the Abelian subalgebra is constituted by  $(n-m)/2$ generators without $I_k$ and $m$
generators containing
$I_k$. If $m< n$ the algebra does not exist for the  real field but, in any case,  the basis can be constructed in such a way that  $c_r^{p,q}=-\bar f^r_{p,q}$. 
In some sense all  these cases could, thus, be considered as Weyl-Drinfeld doubles on $\C$. 

Note that in all the previous expressions for classical Lie algebras the chosen normalization
of the generators of
$s_\pm$ could be changed without destroying the bialgebra structure, since the introduction
of the corresponding inverse factor in $s_\mp$ would be enough. However, in that case the
Weyl-Drinfeld invariance of the double would be broken.

Finally we would like to stress that, in the spirit of Cartan, the
classical series of Lie algebras can be considered as, essentially, the Drinfeld doubles
associated to self-dual Lie bialgebra structures on a very specific set of
solvable Lie algebras (the ones defined by the Cartan generators and the positive roots).

%%%%%%%%%%%%%%% ACKNOWWLEDGMENTS %%%%%%%%%%%
%%%%%%%%%%%%%%%%%%%%%%%%%%%%%%%%%%%%%%%%

\section*{Acknowledgments}

This work was partially supported  by the Ministerio de
Educaci\'on y Ciencia  of Spain (Projects FIS2005-02000 and FIS2004-07913),  by the
Junta de Castilla y Le\'on   (Project VA013C05), and by
INFN-CICyT (Italy-Spain).

%%%%%%%%%%%%%%%%%%%%%%%%%%%%%%%%%%%%%%%%%%%%%%%%%%%%%%%%%%%%%%%%%%%%

%%%%%%%%%%%%%%%%%%%%%% BIBLIOGRAPHY %%%%%%%%%%%%%%%%%%%%%

%%%%%%%%%%%%%%%%%%%%%%%%%%%%%%%%%%%%%%%%%%%%%%%%%%%%%%%%%%%%%%%%%%%%%%%

\end{document}